\documentclass[12pt,reqno]{article}

\usepackage[usenames]{color}
\usepackage{amssymb}
\usepackage{amsmath}
\usepackage{amsthm}
\usepackage{amsfonts}
\usepackage{amscd}
\usepackage{graphicx}

\usepackage[colorlinks=false,
linkcolor=webgreen,
filecolor=webbrown,
citecolor=webgreen]{hyperref}

\definecolor{webgreen}{rgb}{0,.5,0}
\definecolor{webbrown}{rgb}{.6,0,0}

\usepackage{color}
\usepackage{fullpage}
\usepackage{float}

\usepackage{psfig}
\usepackage{graphics}
\usepackage{latexsym}
\usepackage{epsf}
\usepackage{breakurl}

\newcommand{\seqnum}[1]{\href{https://oeis.org/#1}{\rm \underline{#1}}}

\def\modd#1 #2{#1\ \mbox{\rm (mod}\ #2\mbox{\rm )}}

\newcommand{\mean}{\mathop{\mathrm{mean}}}
\newcommand{\lcm}{\mathop{\mathrm{lcm}}}
\newcommand{\li}{\mathop{\mathrm{li}}}
\newcommand{\Exp}{\mathop{\mathrm{Exp}}}
\newcommand{\Gumbel}{\mathop{\mathrm{Gumbel}}}
\newcommand{\itpara}[1]{\smallskip\noindent{\it #1}}
\newcommand{\bfpara}[1]{\smallskip\noindent{\bf #1}}
\newcommand{\rmpara}[1]{\smallskip\noindent{\rm #1}}

\def \BB{{\cal B}}
\def \TT{{\cal T}}

\begin{document}

\noindent~~
\vskip 1cm

\theoremstyle{plain}
\newtheorem{theorem}{Theorem}
\newtheorem{corollary}[theorem]{Corollary}
\newtheorem{lemma}[theorem]{Lemma}
\newtheorem{proposition}[theorem]{Proposition}

\theoremstyle{definition}
\newtheorem{definition}[theorem]{Definition}
\newtheorem{example}[theorem]{Example}
\newtheorem{conjecture}[theorem]{Conjecture}

\theoremstyle{remark}
\newtheorem{remark}[theorem]{Remark}

\begin{center}
\vskip 1cm{\LARGE\bf
On the First Occurrences of Gaps \\  \vskip 1mm
Between Primes in a Residue Class
}
\vskip 0.6cm
{\large Alexei Kourbatov\\
JavaScripter.net \\
Redmond, WA 98052\\
 USA \\
 {\tt akourbatov@gmail.com} \\
\ \\
Marek Wolf\\
Faculty of Mathematics and Natural Sciences\\
Cardinal Stefan Wyszy\'nski University \\
Warsaw, PL-01-938 \\
 Poland \\
 {\tt m.wolf@uksw.edu.pl} \\
}

\vspace{8mm}
{\normalsize\it To the memory of Professor Thomas R.~Nicely (1943--2019)}
\vspace{6mm}
\end{center}

\begin{abstract}\noindent
We study the first occurrences of gaps between primes in the arithmetic progression
(P):~$r$, $r+q$, $r+2q$, $r+3q, \ldots, $
where $q$ and $r$ are coprime integers, $q>r\ge1$.
The~growth trend and distribution of the first-occurrence gap sizes
are similar to those of {\it maximal} gaps between primes in (P).
The histograms of first-occurrence gap sizes, after appropriate rescaling, are well
approximated by the Gumbel extreme value distribution. Computations suggest that
first-occurrence gaps are much more numerous than maximal gaps:
there are $O(\log^2 x)$ first-occurrence gaps between primes in (P) below $x$,
while the number of maximal gaps is only $O(\log x)$.
We explore the connection between the asymptotic density of gaps of a given size and
the corresponding generalization of Brun's constant. For the first occurrence of gap $d$ in (P),
we expect the end-of-gap prime $p\asymp\sqrt{d}\exp({\sqrt{d/\varphi(q)}})$
infinitely often. Finally, we study the gap size as a function of its
index in the sequence of first-occurrence gaps.
\end{abstract}

\section{Introduction}
Let $p_n$ be the $n$-th prime number, and consider the difference between successive primes,
called a {\it prime gap}: $p_{n+1}-p_n$.
If a prime gap is {\it larger} than all gaps before it, we call this gap {\it maximal}.
If a prime gap is either {\it larger} or {\it smaller} than (but not equal to) any gap before it, we call this gap a {\it first occurrence}
\cite{Brent1973,Nicely1999,NicelyNyman,OliveiraeSilva2014,YoungPotler}.

Let $G(x)$ be the largest gap between primes $\le x$:
$$G(x) = \max_{p_{n+1}\le x} (p_{n+1} - p_n).
$$
A lot has been conjectured about large prime gaps; much less has been rigorously proved.
Today's best known results on the growth of $G(x)$ as $x\to\infty$ are as follows:
$$
A\,\textstyle{\log x \log\log x \log\log\log\log x \over \log\log\log x} < G(x)  <  B x^{0.525},
$$
with some positive constants $A$ and $B$.
The upper bound is due to Baker, Harman, and Pintz~\cite{BakerHarmanPintz},
and the lower bound to Ford, Green, Konyagin, Maynard, and Tao \cite{FGKMT}.
Computations suggest that the above double inequality holds with $A=B=1$ for all $x\ge153$,
while the lower bound holds with $A=1$ wherever the left-hand side expression is defined;
cf.\,\cite{FunkhouserGoldstonLedoan}.

The notion of prime gap admits a natural generalization.
Given an increasing integer sequence (S), we can consider gaps $p'-p$ between primes $p$ and $p'$
both of which are in (S), under the condition that $p<p'$ and there are no other primes in (S) between $p$ and $p'$.
Maximal gaps and first occurrence gaps between primes in (S) are also defined in a natural way,
similar to the above.

In our recent work \cite{KourbatovWolf2019} we investigated the statistical properties of maximal gaps
in certain increasing subsequences of the prime numbers.
In this paper we
turn our attention to a related, more numerous family of prime gaps:
the  first occurrences of gaps of a particular size.
As a sequence (S) of interest, here we will choose an arithmetic progression (P): $r+nq$, $n=0,1,2,3,\ldots\,$, with $\gcd(q,r)=1$.
It is plain that the first-occurrence gaps include, as a subsequence, all maximal gaps
between primes $p\equiv r$ (mod $q$);\, e.g., 
sequences \seqnum{A268925} and \seqnum{A330853} from the
{\it On-Line Encyclopedia of Integer Sequences} (OEIS) \cite{oeis}
illustrate the case $q=6$, $r=1$.
We will see that maximal gaps and first-occurrence gaps between primes in (P) share certain statistical properties.

If we limit ourselves to studying the maximal and first-occurrence gaps in the sequence of {\it all primes},
we will find that the available data are quite scarce.
As of 2019, we know only 745 first-occurrence gaps between primes below  $2^{64}$ (OEIS \seqnum{A014320});
of these first occurrences, only 80 gaps are maximal \cite{Nicely_gaps}.
This is one of the compelling reasons to study large gaps between primes in residue classes,
i.e., in arithmetic progressions (P): $r+nq$.
If the value of $q$ used as the common difference of progression (P) is ``not too small'', we get {\em plenty of data} to study large prime gaps.
That is, we get many sequences of first-occurrence gaps corresponding to progressions $r+nq$, for {\it the same fixed} $q$,
which enables us to focus on common features of these sequences.
One such feature is the common histogram of (appropriately rescaled) sizes of first-occurrence gaps.
Other interesting features are the {\em average} numbers of maximal gaps and first-occurrence gaps between primes in (P) not exceeding $x$.

\subsection{Notation}

{\small

\noindent
\begin{tabular}{ll}
$p_n$             & the $n$-th prime;\, $\{p_n\} = \{2,3,5,7,11,\ldots\}$                                  \\
$q$, $r$          & coprime integers, $1\le r<q$                                                           \\
(P) = (P)$_{q,r}$ & the arithmetic progression \ $r, \ r+q, \ r+2q, \ r+3q, \ldots$                        \\
$P_f(d;q,r)$      & the prime starting the first occurrence of gap $d$ between primes in (P)               \\
$P'_f(d;q,r)$     & the prime ending the first occurrence of gap $d$ between primes in (P)                 \\
$p_f(d)$          & the prime starting the first occurrence of gap $d$ in the sequence of all primes       \\
$\gcd(m,n)$       & the greatest common divisor of $m$ and $n$                                             \\
$\lcm(m,n)$       & the least common multiple of $m$ and $n$                                               \\
$\lfloor x\rfloor$& the floor function: the greatest integer $\le x$                                       \\
$\lceil x\rceil$  & the ceiling function: the least integer $\ge x$                                        \\
$\varphi(q)$      & Euler's totient function (OEIS \seqnum{A000010})                                       \\
$f(x)\sim h(x)$   & ${f(x)/h(x)}\to1$ \,as $x\to\infty {\large\phantom{1^{1^1}}}$                          \\
$f(x)\asymp h(x)$ &  $A h(x) \le f(x) \le B h(x)$  for some $A,B>0$                                        \\
$f(x)\lesssim h(x)$    &  $f(x)\sim h(x)$ and $f(x)<h(x)$  for large $x$                                   \\
$A\approx B$      &  $A$ is close to $B$ (the $\approx$ relation is often based on heuristics or numerical data)\\
$\Gumbel(x;\alpha,\mu)$& the Gumbel distribution cdf: \ $\Gumbel(x;\alpha,\mu) =
                    e^{-e^{-{x-\mu\over\vphantom{f}\alpha}}}$                                              \\
$\Exp(x;\alpha)$  & the exponential distribution cdf: \ $\Exp(x;\alpha)=1-e^{-x/\alpha}$                   \\
$\alpha$          & the {\em scale parameter} of exponential/Gumbel distributions, as applicable           \\
$\mu$             & the {\em location parameter} ({\em mode}) of the Gumbel distribution                   \\
$\gamma$          & the Euler--Mascheroni constant: \ $\gamma = 0.57721\ldots$ (OEIS \seqnum{A001620})     \\
$\Pi_2$           & the twin prime constant: \ $\Pi_2 = 0.66016\ldots$ (OEIS \seqnum{A005597}; $\Pi_2^{-1}$ is \seqnum{A167864})       \\
$\log x$          & the natural logarithm of $x$                                                           \\
$\li x$           & the logarithmic integral of $x$: \
                    $\displaystyle\li x \,= \int_0^x{\negthinspace}{{\rm d}t\over\log t}
                                        \,= \int_2^x{\negthinspace}{{\rm d}t\over\log t} + 1.04516\ldots$  \\

                  & {\em Prime counting functions:}                          ${\large\phantom{1^{1^1}}}$   \\
$\pi(x)$          & the total number of primes $p_n\le x$                                                  \\
$\pi(x;q,r)$      & the total number of primes $p=r+nq\le x$                                               \\

                  & {\em Gap counting functions:}                            ${\large\phantom{1^{1^1}}}$   \\
$N_{q,r}(x)$      & the number of maximal gaps $G_{q,r}$ in (P) with endpoints $p\le x$                    \\
$N'_{q,r}(x)$     & the number of first-occurrence gaps in (P) with endpoints $p\le x$                     \\
$\tau_{q,r}(d,x)$ & the number of gaps of a given even size $d=p'-p$ between successive                    \\
                  & primes $p,p'\equiv r$ (mod $q$), with $p'\le x$; \ \
                    $\tau_{q,r}(d,x)=0$ \ if $q\nmid d$ or $2\nmid d$.
                                                                                                           \\
                  & {\em Gap measure functions:}                             ${\large\phantom{1^{1^1}}}$   \\
$G(x)$            & the maximal gap between primes $\le x$: \ $G(x) = \max_{p_{n+1}\le x} (p_{n+1} - p_n)$ \\
$G_{q,r}(x)$      & the maximal gap between primes $p=r+nq \le x$                        \\
$d_{q,r}(x)$      & the latest first-occurrence gap between primes $p=r+nq \le x$ (sect.\,\ref{sect-distribution})             \\
$R(n,q,r)$        & size of the $n$-th record (maximal) gap between primes in (P) \,(sect.\,\ref{S7})      \\
$S(n,q,r)$        & size of the $n$-th first-occurrence gap between primes in (P) \,(sect.\,\ref{S7})      \\
CSG ratio         & the Cram\'er--Shanks--Granville ratio for gap $p'-p$: \
                       $\mbox{CSG}={p'-p \over \varphi(q)\log^2 p'}$ \\
$a(q,x)$   & the expected average gap between primes in (P): \ $a(q,x)=x\varphi(q)/\li x$                  \\
$T_0$, $T_f$, $T_m$
                  & trend functions predicting the growth of large gaps (sect.\,\ref{sec-trend})           \\
\end{tabular}
}

\subsection{Preliminaries}

The following theorems and conjectures provide a broader context for our work on large prime gaps.

{\it The prime number theorem}
\,states that the number of primes not exceeding $x$ is asymptotic to the logarithmic integral $\li x$:
\begin{equation}\label{eq-pnt}
\pi(x) \sim \li x \qquad\mbox{ as } x\to\infty.
\end{equation}

{\it The Riemann hypothesis}\, implies that the error term in (\ref{eq-pnt}) is small:
$$
\pi(x)-{\li x} \,=\, O\!\left(x^{1/2+\varepsilon}\right)
\quad\mbox{ for any } \varepsilon>0;
$$
that is, the numbers $\pi(x)$ and $\lfloor \li x \rfloor$ almost agree in the left half of their digits.

Let $q>r\ge1$ be fixed coprime integers, and consider the arithmetic progression
$$\mbox{(P)} ~~~ r, \ r+q, \ r+2q, \ r+3q, \ \ldots \,.$$

{\it Dirichlet's theorem on arithmetic progressions} \cite{Dirichlet} establishes
that there are infinitely many primes in progression (P) whenever $\gcd(q,r)=1$.
{\it The prime number theorem for arithmetic progressions} \cite{deKoninckLuca}
states that the number of primes in progression (P) not exceeding $x$ is
asymptotically equal to $1/\varphi(q)$ of the total number of primes $\le x$, that is,
\begin{equation}\label{eq-pntap}
\pi(x;q,r) \sim {\pi(x) \over \varphi(q)} \sim {\li x \over \varphi(q)} \qquad\mbox{ as } x\to\infty,
\end{equation}
where $\varphi(q)$ is Euler's totient function.
Under the {\it Generalized Riemann Hypothesis} \cite[p.\,253]{deKoninckLuca}
we have
$$\pi(x;q,r)-{\li x \over \varphi(q)} \,=\, O\!\left(x^{1/2+\varepsilon}\right)
\quad\mbox{ for any } \varepsilon>0.
$$

{\it Polignac's conjecture} \cite{polignac1849}
states that every gap of an even size $d=2n$, $n\in{\mathbb N}$, actually occurs between consecutive primes
infinitely often. (In Section~\ref{sectWhen} we generalize Polignac's conjecture to gaps of size $n\cdot\lcm(2,q)$
between primes in progression (P); see Conjecture \ref{genPolignac}.)

Regarding the behavior of the function $G(x)$---the
maximal prime gap up to $x$---Cram\'er conjectured in the 1930s that $G(x)=O(\log^2 x)$ \cite{Cramer}; \
Shanks set forth a stronger hypothesis: $G(x)\sim\log^2 x$ \cite{Shanks1964}.
Recent models \cite{BanksFordTao,Granville} cast doubt on the latter conjecture;
we still consider it plausible in an ``almost all'' sense (for {\it almost all} maximal gaps).

Let $G_{q,r}(x)$ denote the maximal gap between primes $\le x$ in progression (P).

The \itpara{Generalized Cram\'er conjecture} states that almost all maximal gaps $G_{q,r}(x)$ satisfy
\begin{equation}\label{genCramer}
G_{q,r}(x) < \varphi(q) \log^2 x.   \qquad\mbox{ \cite[eq.~(34)]{KourbatovWolf2019} }
\end{equation}

The \itpara{Generalized Shanks conjecture} is similar: almost all maximal gaps $G_{q,r}(x)$ satisfy
\begin{equation}\label{genShanks}
G_{q,r}(x) \sim \varphi(q) \log^2 x.  \qquad\mbox{ \cite[eq.~(35)]{KourbatovWolf2019}}
\end{equation}

\section{The growth trend of large gaps}\label{sec-trend}

In this section we introduce three trend functions
useful in describing the growth and distribution of large gaps between primes in arithmetic progression (P): $r+nq$.
These  trend functions are asymptotically equivalent to $\varphi(q)\log^2 x$ and differ only in lower-order terms.

\begin{itemize}
\item The function $T_0(q,x)$ is called a {\it baseline trend}.
Roughly speaking, for large enough $x$, this trend curve separates typical sizes of maximal gaps from
typical sizes of non-maximal first-occurrence gaps; see sect.\,\ref{subsec-T0}.

\item The function $T_m(q,x)$ will give us an estimate of the most probable maximal gap sizes.
We have $T_m(q,x)\ge T_0(q,x)$;
see sect.\,\ref{subsec-Tm}.

\item The function $T_f(q,x)$ will give us an estimate of the most probable first-occurrence gap sizes.
For large $x$, we have $T_f(q,x)<T_0(q,x)$;
see sect.\,\ref{subsec-Tf}.
\end{itemize}
Motivated by the prime number theorem for arithmetic progressions (\ref{eq-pntap}),
we give the following estimate $a(q,x)$ of the average gap between primes in (P),
to be used in trend functions.

\begin{definition}
The {\em expected average gap} between primes $p\le x$ in progression (P) is
\begin{equation}\label{aqx}
a(q,x) = {x \varphi(q) \over \li x},
\end{equation}
where $\varphi(q)$ is Euler's totient function, and $\li x$ is the logarithmic integral.
\end{definition}

\subsection{The baseline trend $T_0(q,x)$}\label{subsec-T0}

We now define a trend function that will play a central role in this work.

\begin{definition}
The {\em baseline trend\,} $T_0(q,x)$ for large gaps between primes in (P) is
\begin{equation}\label{T0qx}
T_0(q,x) =a(q,x)\,\log{\li x \over a(q,x)},
\end{equation}
where $a(q,x)$ is given by (\ref{aqx}), and $x$ is large enough, so that ${\li x \over a(q,x)}>1$.
\end{definition}

Formula (\ref{T0qx}) is a concise form of
\cite[eq.\,(33)]{KourbatovWolf2019} in which we set a constant term\footnote{
  The choice of $b=\log\varphi(q)$  for use in \cite[eq.\,(33)]{KourbatovWolf2019} reflects
  our heuristic expectation that, for large $x$, maximal gaps between primes $\le x$ in (P)
  should depend primarily on the quantities $\varphi(q)$, $a(q,x)$, and $\pi(x)$.
  Moreover, this choice is supported by extensive numerical results; cf.\,\cite[sect.\,3.1]{KourbatovWolf2019}.
}
$b=\log\varphi(q)$.
Thus, by definition we have
\begin{equation}\label{T0qx2logli}
T_0(q,x) = {x \varphi(q) \over \li x} \left(2 \log{ \li x \over \varphi(q)} - \log{x\over \varphi(q)} \right).
\end{equation}
(Equations (\ref{T0qx}) and (\ref{T0qx2logli}) can also be viewed as a generalization of the trend equation
derived in \cite{Wolf2011, Wolf2014} for maximal prime gaps $G(x)$, that is, for the particular case $q=2$.)

\begin{remark} The baseline trend $T_0$ defined in (\ref{T0qx}) has the following properties: 
\begin{itemize}
\item[(i)] $T_0(q,x)$ is a continuous and smooth function of $x$ for every fixed $q$. 
\item[(ii)] $T_0(q,x)$ does not depend on the choice of $r$ in progression (P). 
\item[(iii)]  $T_0(q,x)$ does not use unknown parameters determined a posteriori. 
\end{itemize}
\end{remark}

By way of comparison,
for {\it exponentially distributed random} (i.i.d.)~intervals between rare events,
with mean interval $\alpha$ and cdf $\Exp(\xi;\alpha)=1-e^{-\xi/\alpha}$,
the most probable maximal interval in $[0,x]$ is about $\alpha \log{x\over \alpha}$.
Note the absence of $\li x$, in contrast to (\ref{T0qx});
cf.\,\cite[sect.\,4.1.1]{Gumbel}, \cite[sect.\,8]{Kourbatov2013}, and \cite[sect.\,3.3]{Kourbatov2016}.

Clearly, as $x\to\infty$ we have
\begin{equation}
a(q,x) \lesssim \varphi(q)\log x,
\end{equation}
\begin{equation}
T_0(q,x) \lesssim \varphi(q)\log^2 x.
\end{equation}

We numerically experimented with arithmetic progressions (P): $r+nq$\, for many different coprime pairs $(q,r)$, up to $q=10^5$.
Invariably, we found that the majority of maximal gaps between primes $\le x$ in (P) have sizes {\it above} the baseline trend curve $T_0(q,x)$;
that is, we usually have $G_{q,r}(x)>T_0(q,x)$.
At the same time, for large enough $x$, a significant proportion of non-maximal first-occurrence gap sizes are {\it below} the curve $T_0(q,x)$.
These findings hold at least for $x\le10^{14}$.

\subsection{The trend $T_m(q,x)$ of maximal gaps}\label{subsec-Tm}

Most probable {\it maximal} gaps $G_{q,r}(x)$ are close to the following trend
\cite[eqs.\,(33),\,(45)]{KourbatovWolf2019}:
\begin{equation}\label{Tmqx}
T_m(q,x) = T_0(q,x) +  {\mathcal E}_m(q,x),
\end{equation}
where the error term ${\mathcal E}_m(q,x)$ is nonnegative and given by this empirical formula:
\begin{equation}\label{Emqx}
{\mathcal E}_m(q,x) \approx  {b_1  \log\varphi(q) \over (\log\log x)^{b_2}}\cdot a(q,x).
\end{equation}

For $q=2$, formula (\ref{Emqx}) gives ${\mathcal E}_m(2,x)=0$.
Indeed, we found that the baseline trend $T_0(2,x)$, even {\it without} an error term,
satisfactorily describes the most probable sizes of maximal prime gaps $G(x)$.
(Another reasonably good choice of the error term in (\ref{Tmqx}) for $q=2$ is
${\mathcal E}_m(2,x)=ca(2,x)$ with a constant $c\in[0,{1\over2}]$.)

For $q>2$, we found that the error term ${\mathcal E}_m(q,x)$ is positive, and it is best to determine
$b_1$ and $b_2$ in (\ref{Emqx}) {\it a posteriori}.
However, for $x\in[10^7,10^{14}]$ and $q\in[10^2,10^5]$, choosing
\begin{equation}\label{b1b2}
b_1=4, \qquad b_2=2.7
\end{equation}
yields a ``good enough'' trend curve for  typical maximal gap sizes.
This choice of $b_1$, $b_2$ works best  when $q$ is prime or semiprime.
See \cite{Kourbatov2016,KourbatovWolf2019} for further details on maximal gaps.

\subsection{The trend $T_f(q,x)$ of first-occurrence gaps}\label{subsec-Tf}

For the trend equation to better describe the most probable {\it first-occurrence} gap sizes,
we add a lower-order correction term to the baseline trend (\ref{T0qx}):
\begin{equation}\label{Tfqx}
T_f(q,x) = T_0(q,x) + {\mathcal E}_f(q,x).  \\
\end{equation}
Based on extensive numerical data, we find that the term ${\mathcal E}_f(q,x)$ is eventually negative
and given by this empirical formula:
\begin{equation}\label{Efqx}
{\mathcal E}_f(q,x) \approx  (c_0 - c_1 \log \log x)\cdot a(q,x),
\end{equation}
where $a(q,x)$ is defined by (\ref{aqx}).

Like with $T_m(q,x)$, it is best to determine the parameters $c_0$ and $c_1$ in (\ref{Efqx}) {\it a posteriori}.
However, in some special cases, we were able to find predictive formulas for $c_0$ and $c_1$ as functions of $q$:

\begin{itemize}
\item
If $q$ is a {\it prime} $\ge5$, the following empirical formulas
\begin{equation}\label{c0c1primeq}
\begin{split}
c_0 &= c_0(q) \approx 2.18 \, (\log\lcm(2,q))^{0.58}, \\
c_1 &= c_1(q) \approx 1.18 \, (\log\lcm(2,q))^{0.364}
\end{split}
\end{equation}
work well for parameters $c_0$ and $c_1$ in the term ${\mathcal E}_f(q,x)$ (\ref{Efqx}),
for $x\in[10^7,10^{14}]$ and $q\in[5,10^5]$.

\item
If $q$ is an {\it even semiprime} $\ge10$, it is plain that we
deal with the same sequences of even first-occurrence gaps as above, so we can use the same formulas (\ref{c0c1primeq}).

\item
For $q=2$ (gaps in the sequence of all primes), the parameters
\begin{equation}
c_0\approx3, \qquad c_1\approx1.58
\end{equation}
are close to optimal in the term ${\mathcal E}_f(2,x)$, at least for $x<2^{64}$.
\end{itemize}

\subsection{Numerical results}

\begin{figure}[H]
  \centering
  \includegraphics[bb=30 0 540 330,width=5.6in,height=4in,keepaspectratio]{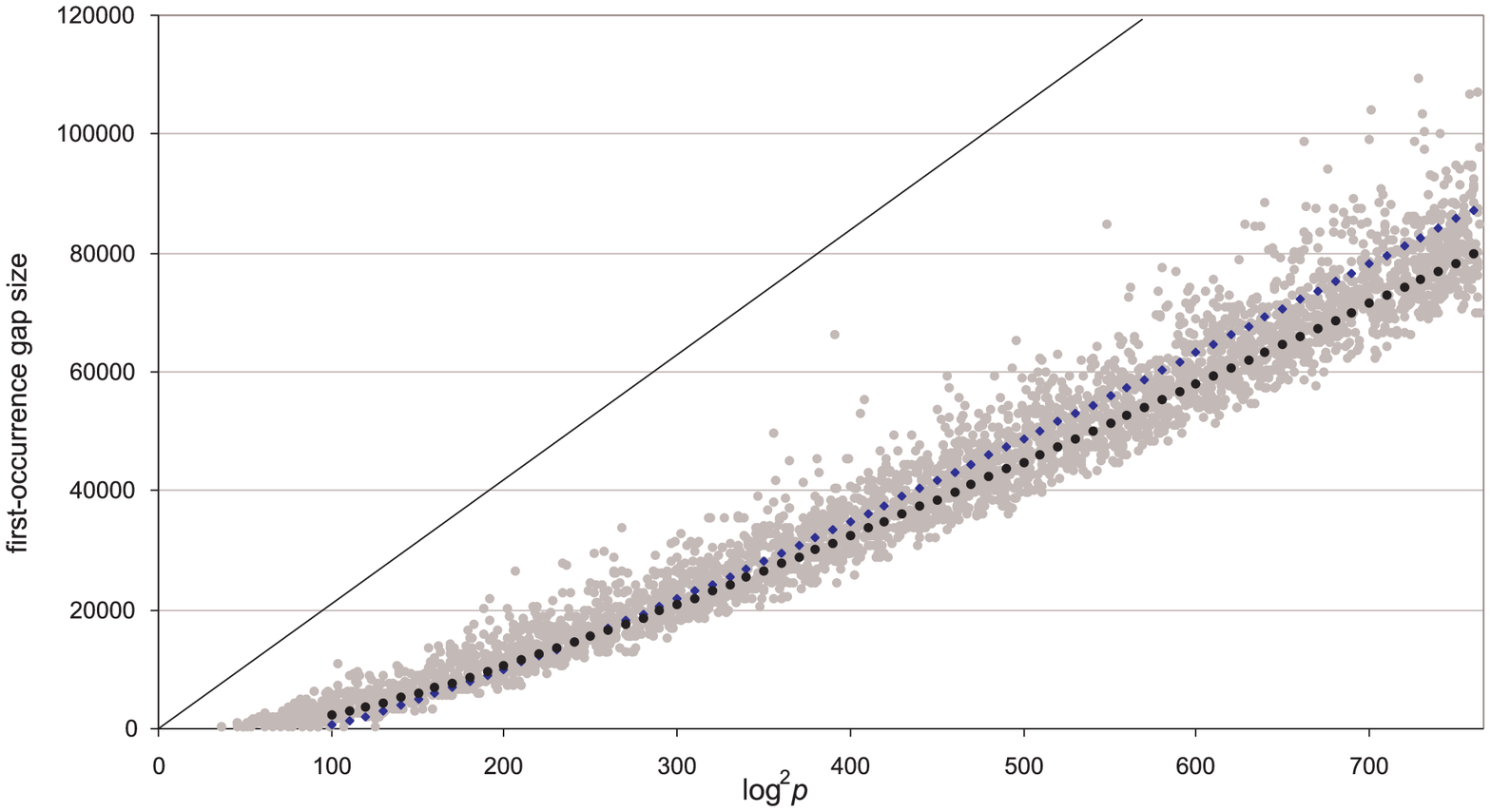}
  \caption{First-occurrence gaps $d_{q,r}$ between primes $p=r+nq\le10^{12}$ for $q=211$, $r\in[1,20]$.
   Dotted curve $(\bullet)$: trend $T_f$ (\ref{Tfqx})--(\ref{c0c1primeq});
 \ blue dotted curve~$(^{_{\color{blue}\,\blacklozenge}})$: baseline trend $T_0$ (\ref{T0qx});
 \ top line: $y = \varphi(q) \log^2 p$.
 The vast majority of {\it maximal} gaps stay above the curve $T_0$, while
 a growing proportion of non-maximal first-occurrence gaps are observed below $T_0$.
}
\label{figq211trend}

\vspace*{\floatsep}

\begin{center}
\includegraphics[bb=-6 0 800 180,width=8.65in,height=2in,keepaspectratio]{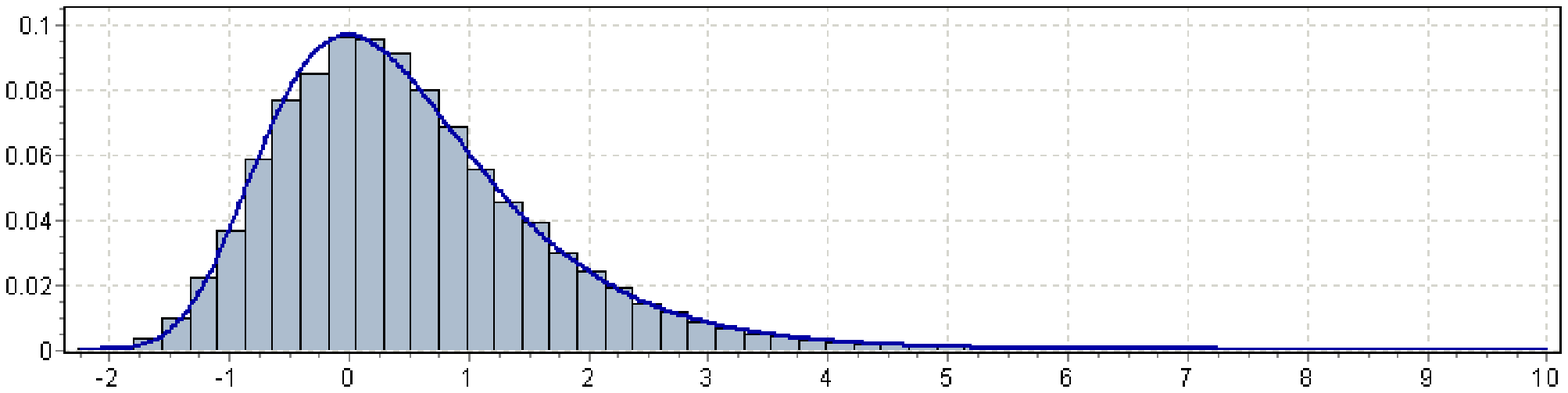}
\caption{Histogram of rescaled values $u$ (\ref{rescaledu}) for the first-occurrence gaps
 between primes $p=r+nq\in[10^7,10^{12}]$ for $q=211$ and all $r\in[1,210]$.
 Bins:~53,\, bin size $\approx0.23$.
 Smooth curve: best-fit Gumbel distribution (pdf) with scale $\alpha=0.87692$ and mode $\mu=-0.00902$.
}
\label{figq211dist}
\end{center}
\end{figure}

Using {\tt PARI/GP} \cite{Pari} we performed computational experiments with progressions (P): $r+nq\,$
for many different values of $q\in[2,10^5]$.
Figure \ref{figq211trend} shows the sizes of first-occurrence gaps between primes in twenty progressions (P),
for $q=211$ and $r\in[1,20]$.
Results for other values of $q$ look similar to Fig.\,\ref{figq211trend}.
All our numerical results show that trend curves of the form (\ref{Tfqx}), (\ref{Efqx}) satisfactorily describe the growth of
typical first-occurrence gaps between primes in progressions (P).

$$~$$

\section{The distribution of first-occurrence gap sizes}\label{sect-distribution}

Extensive numerical data allows us to observe that typically
a first-occurrence gap is also the largest gap between primes in the interval $[x,x+x/(\log x)^{1+\varepsilon}]$,
for some $\varepsilon>0$. This invites the question (cf.\,\cite{Kourbatov2013}--\cite{LiPrattShakan2017}):
do first-occurrence gap sizes, after an appropriate rescaling transformation, obey
the same distribution as maximal gap sizes---the Gumbel extreme value distribution?

Let $d=d_{q,r}(x)$ denote the size of the first-occurrence gap between primes $p\in[1,x]\cap$(P);
specifically, the gap $d_{q,r}(x)$ refers to the first occurrence
{\it closest to} $x$ on the interval $[1,x]$.
Infinitely many distinct values of $d_{q,r}(x)$ are maximal gaps, for which $d_{q,r}(x)=G_{q,r}(x)$.
At the same time, infinitely many distinct values of $d_{q,r}(x)$ are {\it not} maximal gaps, and
for these first occurrences we have $d_{q,r}(x)<G_{q,r}(x)$.

\subsection{Rescaling transformation}

Extreme value theory motivates the following rescaling transformation (cf.\,\cite[sect.\,3.2]{KourbatovWolf2019}):
\begin{equation}\label{rescaledu}
\mbox{ gap size } d_{q,r}(x) ~~\mapsto~~ u = {d_{q,r}(x) - T_f(q,x) \over a(q,x)},
\end{equation}
where
$a(q,x)$ is the expected average gap (\ref{aqx}), \
$T_f(q,x)$ is the trend function of first-occurrence gaps (\ref{Tfqx}),
and $x$ runs through all end-of-gap primes in a given table of first-occurrence gaps.
Every distinct gap size $d_{q,r}(x)$ is counted only once and mapped to its rescaled value $u$.

Figure \ref{figq211dist} shows a combined histogram of rescaled values $u$ for first-occurrence gaps between
primes in progressions (P) for $q = 211$ and all $r\in[1,210]$.
We can see at once that the histogram is skewed to the right.
A good fit is the Gumbel extreme value distribution---which typically occurs in statistics as the distribution of maxima
of i.i.d.\ random variables with an initial {\it unbounded light-tailed} distribution.

\subsection{\bf Distribution tails}
Observe that the {\it left tail} of the histogram of rescaled gaps is light and short.
That is, first-occurrence gaps are {\it never too far below their trend} $T_f(q,x)$;
recall that the trend function $T_f(q,x)$ itself is asymptotic to $\varphi(q)\log^2 x$.
This leads us to the following

\begin{conjecture}\label{liminf}
Let $d_{q,r}(x)$ be the first-occurrence gap function defined above.
Then
$$
\liminf\limits_{x\to\infty} {d_{q,r}(x) \over \varphi(q)\log^2 x} =1.
$$
\end{conjecture}

On the other hand, the {\it right tail} of the rescaled gaps histogram is long. For the sequence of all primes,
Cram\'er's model \cite{Cramer} suggests $\limsup{d_{2,1}(x) \over \log^2 x} =1$.
Banks, Ford, and Tao \cite{BanksFordTao} discuss models yielding a larger upper limit,
$2e^{-\gamma}\approx1.1229$;  \,cf.~Granville \cite[p.\,24]{Granville}.

\subsection{Goodness of fit}
For best-fit Gumbel distributions, we found that the Kolmogorov--Smirnov (KS) goodness-of-fit statistic
was usually not much better (and sometimes worse)
than the corresponding KS statistic measured for the data subset including {\it maximal gaps only} \cite{KourbatovWolf2019}.
For example, the KS statistic for the histogram shown in Fig.\,\ref{figq211dist} is 0.00815
(vs.~0.00628 for the corresponding data subset comprising maximal gaps only).

\subsection{Distribution scale}
For data sets containing all first occurrence gaps---as well as for maximal gaps only---the scale parameter $\alpha$ of best-fit Gumbel distributions
was in the range $\alpha\in[0.7,1]$. In both cases, when we moved from smaller to larger end-of-gap primes,
the scale $\alpha$ of best-fit Gumbel distributions appeared to slowly grow towards 1.

\subsection{Beyond Gumbel}
Among three-parameter distributions, a very good fit for first-occurrence gap histograms is the
{\it Generalized Extreme Value} (GEV) distribution.
Typical best-fit GEV distributions for rescaled values (\ref{rescaledu}) have a small negative {\it shape parameter}, usually in the interval $[-0.04,0]$.
Recall that the Gumbel distribution is a GEV distribution whose shape parameter is exactly zero.  From extreme value theory
we also know that the GEV distribution of maxima with a negative shape parameter occurs
for i.i.d.\ random variables with an initial distribution {\it bounded from above}.

On the other hand, if we fitted a GEV distribution to a data subset including {\it maximal gaps only}, then
the resulting shape parameter was very close to zero,
which is typical for random variables with an {\it unbounded light-tailed} distribution (the Gumbel type).

\subsection{Is there a limiting distribution?}
Our observations are compatible with the existence of a Gumbel limit law for
(appropriately rescaled) maximal gaps, as well as for first-occurrence gaps between primes in progressions (P).
However, the convergence to the hypothetical Gumbel limit law appears to be better for data sets comprising maximal gaps only,
as opposed to data sets including all first-occurrence gaps (in a given range of primes in progressions (P) for a given $q$).

An alternative hypothesis is that the Gumbel limiting distribution exists for rescaled maximal gaps, but not for first-occurrence gaps.
Of course it is also possible that there is no limiting distribution at all, and the Gumbel distribution is simply a good approximation
of the rescaled gap values.

\section{When do we expect the first occurrence of gap $d$?}\label{sectWhen}

Let $d$ be the size of a gap between primes in progression (P): $r+nq$, and denote by
$$
P_f(d;q,r) \quad\mbox{ and }\quad P'_f(d;q,r)
$$
the primes that, respectively, start and end the {\it first occurrence} of gap $d$.
So we have $$d = P'_f(d;q,r)-P_f(d;q,r).$$

\begin{conjecture}\label{genPolignac}
{\it Generalization of Polignac's conjecture for an arithmetic progression} (P):
$r$,~$r+q$, $r+2q,\ldots$\,\,.
For every pair of coprime integers $(q,r)$
and every $d=n\cdot\lcm(2,q)$, $n\in{\mathbb N}$, we have $$P_f(d;q,r) < \infty.$$
Moreover, every gap of size $n\cdot\lcm(2,q)$, $n\in{\mathbb N}$, occurs {\it infinitely often}
between primes in progression (P).
\end{conjecture}

\begin{conjecture}
Let $d$ be a gap between primes in (P). There exists a constant $C$ such that
for the sequence of first-occurrence gaps $d$ infinitely often we have
\begin{equation}\label{dblineq-Pfdqr}
e^{\sqrt{d/\varphi(q)}} < P_f(d;q,r) < C\sqrt{d}\,e^{\sqrt{d/\varphi(q)}},
\end{equation}
and the difference $P_f(d;q,r) - C\sqrt{d}\,e^{\sqrt{d/\varphi(q)}}$
changes its sign infinitely often as $d\to\infty$.
\end{conjecture}

For maximal gaps between primes in (P), we state the following {\it stronger} conjectures.

\begin{conjecture}
Let $d=G_{q,r}(p')$ be a {\it maximal gap} between primes in (P).
Then there exists a constant $C$ such that a {\it positive proportion} of maximal gaps satisfy
the double inequality~(\ref{dblineq-Pfdqr}),
and the difference $P_f(d;q,r) - C\sqrt{d}\,e^{\sqrt{d/\varphi(q)}}$ changes its sign infinitely often as $d\to\infty$.
\end{conjecture}

\begin{conjecture}
Let $d=G_{q,r}(p')$ be a {\it maximal gap} between primes in (P).
Then, as $d\to\infty$, a {\it positive proportion} of maximal gaps satisfy
\begin{equation}\label{Pfdqr-asymp}
P_f(d;q,r) \asymp \sqrt{d}\,e^{\sqrt{d/\varphi(q)}}.
\end{equation}
In other words, there are constants $C_0,C_1>0$ such that a positive proportion of maximal gaps satisfy
the double inequality
\begin{equation}\label{dblineq-C0C1}
C_0\sqrt{d}\,e^{\sqrt{d/\varphi(q)}} < P_f(d;q,r) < C_1\sqrt{d}\,e^{\sqrt{d/\varphi(q)}}.
\end{equation}
\end{conjecture}

For primes up to $10^{14}$ and $q\le10^5$, computations suggest that if we take $C_0={1\over10}$, $C_1=10$,
then the majority of maximal gap sizes will satisfy (\ref{dblineq-C0C1});
a sizable proportion of non-maximal first-occurrence gaps will also be found within the bounds (\ref{dblineq-C0C1}).

How did we arrive at the estimates (\ref{dblineq-Pfdqr})--(\ref{dblineq-C0C1})?
The left-hand side $e^{\sqrt{d/\varphi(q)}}$ in (\ref{dblineq-Pfdqr}) stems from the conjecture
\cite[eq.\,(34)]{KourbatovWolf2019}
that {\it almost all} maximal gaps between primes in (P) satisfy
$$G_{q,r}(p) < \varphi(q) \log^2p.$$
(This is our generalization of Cram\'er's conjecture \cite{Cramer}.)

To see why we have $C\sqrt{d}\,e^{\sqrt{d/\varphi(q)}}$ in the right-hand side of (\ref{dblineq-Pfdqr})--(\ref{dblineq-C0C1}), consider the function
\begin{equation}\label{tildePdq}
\tilde{P}(d,q) = \sqrt{d}\,e^{\sqrt{d/\varphi(q)}}.
\end{equation}
For any fixed $q\ge2$, one can check that
\begin{equation}\label{limittildePdq}
\lim_{x\to\infty}{\tilde{P}(T_0(q,x),q) \over x} = e^{-1/2} \approx 0.60653.
\end{equation}
That is, for large $d$ and $x$, the function $e^{1/2}\tilde{P}(d,q)$ is ``almost inverse'' of
$T_0(q,x)$ in (\ref{T0qx}).

\begin{remark}
We do not know whether the lower bound in \eqref{dblineq-Pfdqr} is reversed infinitely often or finitely often.
For the special case $q=2$ (gaps in the sequence of all primes),
heuristic models \cite{BanksFordTao,Granville} lead to
$\limsup\limits_{x\to\infty}{G(x)\over\log^2 x}\ge 
2e^{-\gamma}\approx1.1229$---suggesting that, for first-occurrence prime gaps, the lower bound in \eqref{dblineq-Pfdqr} is reversed infinitely often.
However, it is quite difficult to find examples of progressions (P): $r + nq$ with exceptionally early first occurrences
of gaps $d$ such that $P'_f(d;q,r)<e^{\sqrt{d/\varphi(q)}}$ and the Cram\'er--Shanks--Granville ratio
$\mbox{CSG}={d_{q,r}(x)\over\varphi(q)\log^2 x}>1$.
We know only a modest number of maximal gaps of this kind
\cite[Table 1]{KourbatovWolf2019}, \cite{Raab2020}--\cite{Raab2020may}.
We have never seen examples of sign reversal of the lower bound in \eqref{dblineq-Pfdqr}
 ($\mbox{CSG}>1$) for {\it non-maximal} first-occurrence gaps in (P).
Note also that Sun \cite[conj.\,2.3]{Sun2013} proposed a hypothesis similar to
Firoozbakht's conjecture for gaps between primes in (P),
which would imply (at most) a finite number of such sign reversals for any given progression (P).
\end{remark}

\section{A detour: generalization of Brun's constants}\label{SectBrunConst}

We will now give another (less direct) heuristic way of deriving
the estimate $C\sqrt{d}\,e^{\sqrt{d/\varphi(q)}}$ in (\ref{dblineq-Pfdqr})--(\ref{dblineq-C0C1}).
Our approach is to estimate the sum of reciprocals of prime pairs $(p,\,p+d)$
separated by a given fixed gap $d$.
(Wolf \cite{Wolf1997} used this heuristic approach to treat the case $q=2$,
 obtaining the estimate\footnote{
 The preprint \cite{Wolf1997} employs the ``physical'' notation where the relation $f\sim g$
 is a shorthand for ``the functions $f$ and $g$ have the same order of magnitude.''
 Thus $\sqrt{d}\,e^{\sqrt{d}}$ is just an order-of-magnitude estimate.
}
$\sqrt{d}\,e^{\sqrt{d}}$ for the first occurrence of gap $d$.)

First, let us recall some definitions.
{\it Twin primes} are pairs of prime numbers $(p,\,p+2)$ separated by gap 2---the smallest possible gap between odd primes.
{\it Cousin primes} are pairs of prime numbers $(p,\,p+4)$ separated by gap 4.

Brun \cite{Brun} proved in 1919 that the series consisting of reciprocal twin primes converges.
{\it Brun's constant} $B_2$ for twin primes is the sum of reciprocals
\cite[\seqnum{A065421}]{Fry_Nesheiwat_Szymanski,Nicely_Brun,Nicely2001,SebahGourdon}
$$
B_2 =
\left({1\over3}+{1\over5}\right) +
\left({1\over5}+{1\over7}\right) +
\left({1\over11}+{1\over13}\right) + \cdots \approx 1.90216.
$$
A generalization of Brun's theorem exists for series of reciprocal values of prime pairs $(p,\,p+d)$ separated by {\it any fixed gap} $d$ \cite{Segal}.
For example, Brun's constant $B_4$ for cousin primes (with the pair $(3,7)$ omitted) is
\cite[\seqnum{A194098}]{Wolf1996}
$$
B_4 =
\left({1\over7}+{1\over11}\right) +
\left({1\over13}+{1\over17}\right) +
\left({1\over19}+{1\over23}\right) + \cdots \approx 1.19704.
$$

Note that twin primes $p$, $p+2$ are always  consecutive.
Except for the pair $(3,7)$, primes $p$, $p+4$ are also necessarily consecutive.
If $d\ge6$, however, there may also be other primes between the primes $p$ and $p+d$.
In what follows, we will consider gaps $d\ge6$.

Suppose that primes $p$ and $p'$ are in the residue class $r$ (mod $q$), $p<p'$,
and there are {\it no other primes} in this residue class between $p$ and $p'$.
If $ d=n\cdot\lcm(2,q)$, then
the number of gaps of size $d=p'-p$ between primes $p,p'\equiv r$ (mod $q$), with $p'\le x$, is roughly
\begin{equation}\label{main2}
\tau_{q,r}(d,x) \approx C_2\prod_{p|d\atop p>2}\frac{p-1}{p-2} \cdot \frac{ \pi^2(x)}{\varphi ^2(q)x} \, e^{-{d\cdot\pi(x)\over\varphi(q)x}},
\end{equation}
cf.\,\cite[eq.~(2)]{GoldstonLedoan2011} and \cite[eq.~(30)]{KourbatovWolf2019}.
In formula (\ref{main2}), $p$ runs through odd prime factors of $d$,
\begin{equation}\label{params-c-s}
C_2={c\over s}, \qquad
c=\lcm(2,q), \qquad
s= \mean\limits_{d\,=\,nc \atop n\in\mathbb N} \, \prod_{p|d\atop p>2}\frac{p-1}{p-2}
= \Pi_2^{-1} \, \prod_{p|q\atop p>2}\frac{p}{p-1},
\end{equation}
and $\Pi_2={\prod\limits_{p>2}{p(p-2)\over(p-1)^2}}\approx0.66016$ is the twin prime constant.
In particular,  $s=\Pi_2^{-1}\approx1.51478$ if $q$ is a power of 2; see {\it Appendix}.
Note that
\begin{equation}\label{tauzero}
\tau_{q,r}(d,x) = 0 \quad\mbox{ if } 2\nmid d \mbox{ ~or~ } q\nmid d.
\end{equation}

Let ${\cal T}_{d;q,r}$ denote the set of primes $p,p'$ at both ends of the gaps $d$ counted in (\ref{main2}):
\begin{equation}
{\cal T}_{d;q,r} \ =\  \{p, p': ~   p' - p=d  ~\mbox{ and }~ p, p'  \equiv r \!\!\!\! \mod q\}.
\end{equation}
Consider the sum of reciprocals of all primes in ${\cal T}_{d;q,r}$:
\begin{equation}
B_d(q,r)\,=\sum_{p\,\in {\cal T}_{d;q,r}}  \frac{1}{p}.
\label{def_Bd}
\end{equation}
We adopt the rule that if a given gap $d$ appears twice in a row:
$d=p'-p=p''-p'$, then the corresponding middle prime $p'$ is counted twice.
The real numbers $B_d(q,r)$ defined by (\ref{def_Bd}) can be called {\it the generalized Brun constants}
for residue class $r$ (mod $q$).

Let us define the partial (finite) sums:
\begin{equation}\label{Bdxqr}
B_d(x;q,r)\,=\sum_{p\in \TT_{d;q,r}\atop p\le x} \frac{1}{p}.
\end{equation}
Hereafter we will work with rough estimates of partial sums (\ref{Bdxqr}).
To distinguish the actual partial sum $B_d(x;q,r)$ from its estimated value,
we will write $\BB_d(x;q)$ for the estimate.
Likewise, for a rough estimate of the constant (\ref{def_Bd}),
we will write $\BB_d(\infty;q)$.
Our estimates ${\cal B}_d(\,\cdot\,;q)$ will be independent of $r$.

Given a gap size $d=p'-p$, observe that the density of pairs $(p,p')\subset\TT_{d;q,r}$ near $x$ is roughly $\tau_{q,r}(d,x)/x$.
Putting $\pi(x)\approx x/\log(x)$ in equation (\ref{main2}), we obtain
\begin{equation}\label{integralB}
\BB_d(x;q)
\approx
\sum_{p\in \TT_{d;q,r}} \frac{1}{p} ~- \sum_{p\in \TT_{d;q,r}\atop p>x} \frac{1}{p}
\approx
\BB_d(\infty,q) - 2C_2  \prod_{p\mid d} \frac{p-1}{p-2}\int_x^\infty \!\!\frac{e^{-{d\over\varphi(q)\log u}}}{\varphi^2(q)\,u\log^2 u}\, {\rm d}u.
\end{equation}
After the substitution $v={d \over \varphi(q)\log u}$ the integral (\ref{integralB}) can be calculated explicitly:
\begin{equation}
\BB_d(x;q) \approx \BB_d(\infty;q) -
\frac{2C_2}{\varphi(q)d} \prod_{p\mid d} \frac{p - 1}{p - 2} \left(1-e^{-{d\over\varphi(q)\log x}}\right).
\label{Bd_x}
\end{equation}
We will require that
$\BB_d(x;q)\to0$ as $x\to1^+$ (indeed, the actual partial sum $B_d(x;q,r)$ will be zero up to
the first occurrence of the gap $d$). Taking the limit $x\to 1^+$ in (\ref{Bd_x}) we obtain
\begin{equation}
\BB_d(\infty;q) \approx \frac{2C_2}{\varphi(q) d} \prod_{p\mid d} \frac{p - 1}{p - 2}.
\label{main-B-d}
\end{equation}
Then the estimate $\BB_d(x;q)$, as a function of $x$ for fixed $q$ and $d$, has the form
\begin{equation}
\BB_d(x;q) \approx  \frac{2C_2}{\varphi(q)d}  \prod_{p\mid d} \frac{p - 1}{p - 2} \,e^{-{d\over\varphi(q)\log x}}.
\label{Bd-od-x}
\end{equation}
Recall that, for $d=n\cdot\lcm(2,q)$, the mean value of
$\prod_{p\mid d} \frac{p - 1}{p - 2}$ is $s$;
see (\ref{main2}), (\ref{params-c-s}).
Therefore, we can skip the ratio $\prod_{p\mid d} \frac{p - 1}{p - 2}/s$ (which is near 1 infinitely often) and get
\begin{equation}
\BB_d(x;q) \approx \frac{2c}{\varphi(q)d} \, e^{-{d\over\varphi(q)\log x}} =
\frac{A}{d} \, e^{-{d\over\varphi(q)\log x}},
\end{equation}
where $c=\lcm(2,q)$ and $A=2c/\varphi(q)=O(\log\log q)$; moreover\footnote{
 We easily check that $A=4$ when $q=2$, and $A\le6$ when $q$ is an odd prime.
 On the other hand, one can also show that $A=O(\log\log q)$ when $q$ is a primorial (which constitutes an extreme case;
 see \cite{Nicolas1983}).
}
$A=O(1)$ for prime $q$.

Suppose that gap $d$ occurs for the first time at $\xi=P_f(d;q,r)$.
We can estimate $\xi$ from the condition $\BB_d(\xi;q)\approx2/\xi$, so we have the following equation for $\xi$:
\begin{equation}\label{Aexpd2xi}
\frac{A}{d} \, e^{-{d\over\varphi(q)\log\xi}} = \frac{2}{\xi}.
\end{equation}
Taking the log of both sides of (\ref{Aexpd2xi}) and discarding the small term $\log{A\over2}$, we get a quadratic equation for $t=\log\xi$:
\begin{equation}\label{quadratic}
t^2 - t\log d - {d/\varphi(q)} = 0.
\end{equation}
Its positive solution is $t\approx{1\over2}\log d+\sqrt{d/\varphi(q)}$,
which allows us to conjecture that
\begin{equation}\label{dexpsqrtdphi}
\xi = P_f(d;q,r) \asymp  \sqrt{d} \, e^{\sqrt{d/\varphi(q)}}
\end{equation}
infinitely often. 

\begin{remark}
Setting $q=2$ in (\ref{dexpsqrtdphi}), for the first occurrence of gap $d$ in the sequence of all primes we get
$p_f(d) \asymp \sqrt{d} \, e^{\sqrt{d}}$; cf.\,\cite{Wolf1997}.
Thus, in  progressions (P) with large $q$, gaps of a given size $d$ appear much earlier
(due to the division by $ \varphi(q)$ under square root in the exponent).
Note, however, that the ``natural unit'' of gaps in progressions (P) is larger than
that for all primes; namely, all even gaps between primes in (P) have sizes $n\cdot\lcm(2,q)$, $n\in{\mathbb N}$.
\end{remark}

\section{How many first-occurrence gaps are there below $x$?}

\begin{figure}[H]

\begin{center}
  \includegraphics[width=5in,height=4.6in,keepaspectratio]{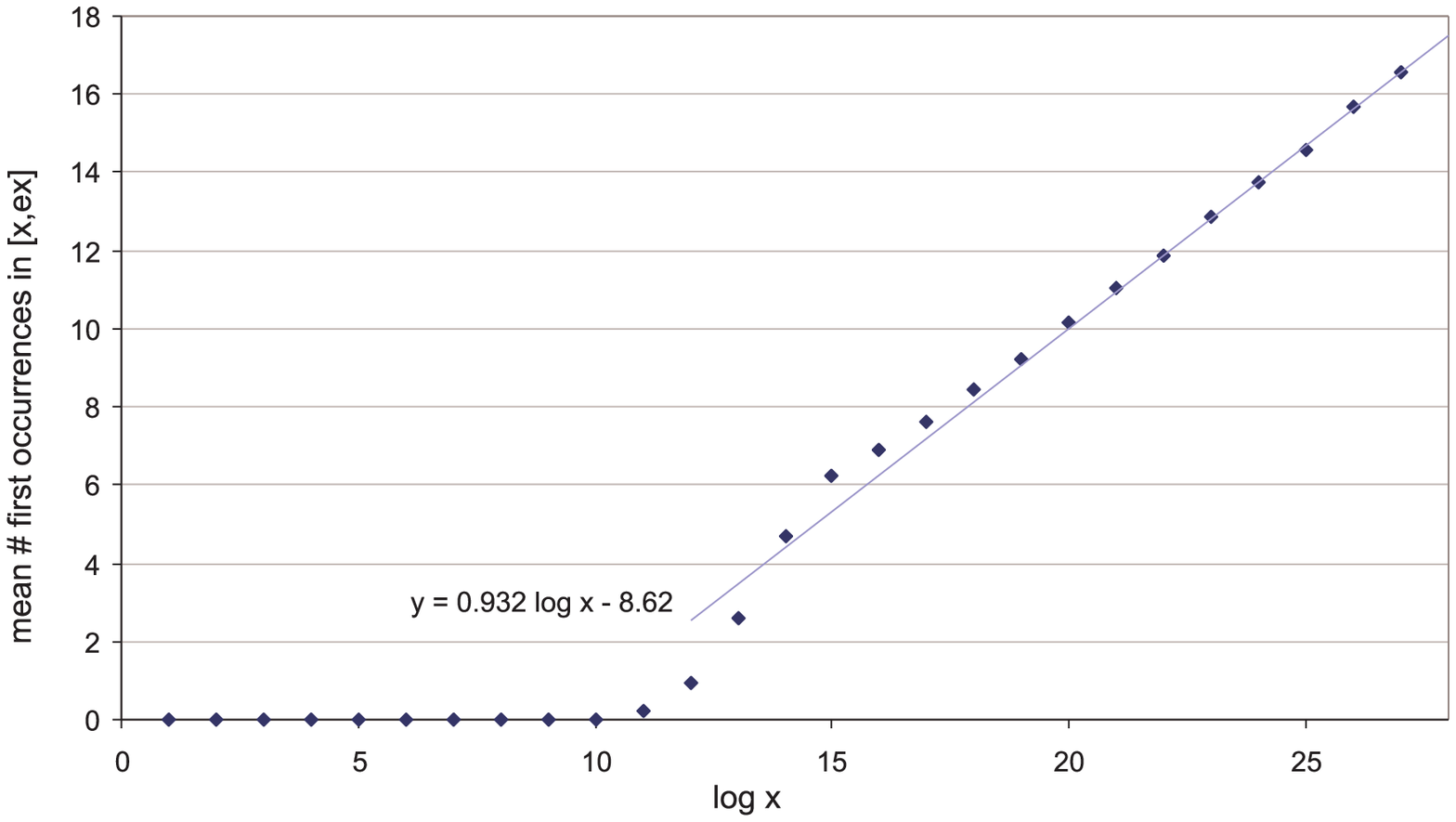}
  \caption{Mean number of first-occurrence gaps between primes $p=r+nq\in[x,ex]$,\,
   for $q=17011$, $x=e^j$, $j\le27$.
  Averaging for all $r\in[1,17010]$.
  }
  \label{figq17011meanfog}

\vspace*{8mm}

\includegraphics[width=5in,height=2.94in,keepaspectratio]{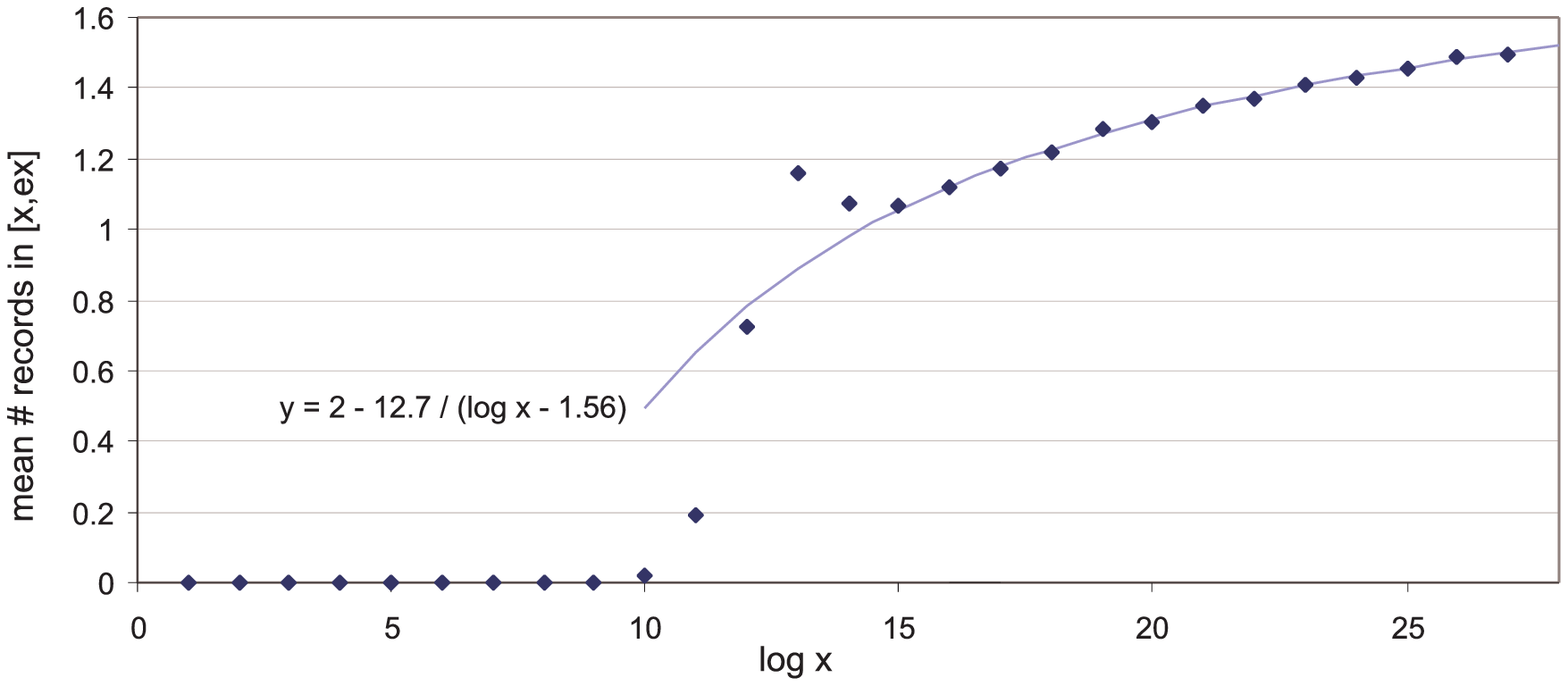}
\caption{Mean number of record (maximal) gaps between primes $p=r+nq\in[x,ex]$,\,
  for $q=17011$, $x=e^j$, $j\le27$.
  Averaging for all $r\in[1,17010]$.
}
\label{figq17011meanrec}
\end{center}
\end{figure}

Using a modified version of our {\tt PARI/GP} code from \cite{KourbatovWolf2019}
we computed the mean numbers of first-occurrence gaps between primes $p=r+nq$
in intervals $[x,ex]$, $x=e^j$, $j=1,2,\ldots,27$.
Figure \ref{figq17011meanfog} shows the results of this computational experiment for $q=17011$.
For comparison, we also computed the corresponding mean numbers of record (maximal) gaps
between primes $p=r+nq$ in the same set of intervals $[x,ex]$; the results are shown in Fig.\,\ref{figq17011meanrec}.
The computation resulting in Figs.\,\ref{figq17011meanfog}, \ref{figq17011meanrec} took about 30 hours
on an Intel Core i7-8550U CPU.

We found that, for large enough $x$, the mean number of first-occurrence gaps between primes $r+nq\in[x,ex]$
grows like a linear function of $\log x$; see Fig.\,\ref{figq17011meanfog}.
It is reasonable to expect that the {\it total} number $N'_{q,r}(x)$ of first-occurrence gaps
between primes $\le x$ in progression (P): $r+nq$ can be approximated by a quadratic function of $\log x$, and for large $x$ we will have
\begin{equation}
N'_{q,r}(x)
\approx {T_0(q,x)\over\lcm(2,q)}
\lesssim {\varphi(q)\log^2x\over\lcm(2,q)}
= O(\log^2x),
\end{equation}
where $T_0(q,x)$ is given by (\ref{T0qx}).

Thus the number $N'_{q,r}(x)$ of first-occurrence gaps between primes in (P) is asymptotically {\it much greater} than the number
$N_{q,r}(x)$ of record (maximal) gaps in (P).
Indeed, the latter seems to be only $O(\log x)$.
More precisely, our computations (see Fig.\,\ref{figq17011meanrec}) suggest that
the average number of maximal gaps observed in $[x,ex]$ is
\begin{equation}\label{eq-hyperbola}
\mean\limits_{r} \big( N_{q,r}(ex) - N_{q,r}(x) \big)
 \approx\, 2 - {\kappa(q) \over \log x + \delta(q)}.
\end{equation}
Accordingly, we conjecture that the number of maximal gaps between primes $\le x$ in (P) is
\begin{equation}\label{eq-2logx}
N_{q,r}(x) \sim 2\log x = O(\log x) \qquad\mbox{ as }x\to\infty.
\end{equation}
The constant 2 in (\ref{eq-hyperbola}) and (\ref{eq-2logx}) can be justified by a heuristic argument
based on two hypotheses: the generalized Shanks conjecture (\ref{genShanks}) and
the assumption that, on average, consecutive record gaps near $x$ differ by
the {\it average} gap $\sim \varphi(q)\log x$ between primes in (P).
(See \cite{Kourbatov2018, KourbatovWolf2019} for details.)

\section{Maximal gaps and first occurrences as functions of their sequential number}\label{S7}

As before, consider the arithmetic progression (P): \ $r, \ r+q, \ r+2q, \ r+3q,\,\ldots$~.

\begin{conjecture}\label{conj-secfopgap}
Let $S(n;q,r)$ be the size of the $n$-th first-occurrence gap between primes in progression (P).
Then the sizes $S(n;q,r)$ of almost all first-occurrence gaps satisfy
\begin{equation}\label{dblSnqr}
n \le S(n;q,r) \le 2nq\lceil \log^2 q \rceil.
\end{equation}
Moreover, all gap sizes divisible by $\lcm(2,q)$ occur early enough, so that
\begin{equation}
S(n;q,r) \approx n\cdot\lcm(2,q)  = O_q(n) \qquad\mbox{ as }n\to\infty.
\end{equation}
{\it Example}: for $q=6$ and $r=5$, we have $S(n;6,5)=$\,\seqnum{A334543}$(n)$.
\end{conjecture}

For comparison, here is a similar conjecture for {\it maximal} gap sizes:

\begin{conjecture}\label{conj-secmaxgap}
Let $R(n;q,r)$ be the size of the $n$-th record (maximal) gap between primes in (P).
Then the sizes $R(n;q,r)$ of almost all maximal gaps satisfy \cite{Kourbatov2018,puzzle_84}
\begin{equation}\label{dblRnqr}
{\varphi(q) n^2 \over 6} <  R(n;q,r) < \varphi(q) n^2 + (n+2)q\log^2 q.
\end{equation}
\end{conjecture}
We do not know any exceptions to (\ref{dblRnqr}); in particular, there are no exceptions for $q\le2000$, $n<15$.
Heuristic reasoning \cite{Kourbatov2018, KourbatovWolf2019} suggests that for large $n$ we will have
\begin{equation}
R(n;q,r) \approx {\varphi(q) n^2 \over 4} = O_q(n^2) \qquad\mbox{ as }n\to\infty.
\end{equation}
(At the same time, for small $n$, the ratio $R(n;q,r)\over\varphi(q) n^2$ may be greater than 1.)

\itpara{Example}: for $q=6$ and $r=5$, we have $R(n;6,5)=$\,\seqnum{A268928}$(n)$.

\medskip
Conjectures \ref{conj-secfopgap} and \ref{conj-secmaxgap} are another way to restate our observation
that first-occurrence gaps between primes in progression (P) are much more numerous than maximal gaps.

\begin{remark}
Setting $n=0$ in (\ref{dblRnqr}), we get that the ``0th record gap'' is less than $2q \log^2 q$,
which can be interpreted as an almost-sure upper bound on the least prime in progression (P); cf.\,\cite[sect.2]{LiPrattShakan2017}.
That is, the ``0th record gap'' between primes in (P) is the distance from 0 to the least prime in (P),
which is (almost always) bounded by $2q \log^2 q$.
\end{remark}

\section{Appendix}
In this appendix we calculate the average value of the product
\begin{equation}\label{defSd}
S(d) = \prod_{p|d\atop p>2} {p-1\over p-2}
\end{equation}
when $d$ runs through a given arithmetic progression, $d=r+nq$, $n\in\mathbb N$. So our goal is to find
$$
 \mean\limits_{d=r+nq\atop n\in\mathbb N} S(d) = \lim_{N\to\infty} {1\over N} \sum_{n=1}^N S(r+nq).
$$
For generality, in the appendix, $q$ and $r$ are not necessarily coprime, and $r$ is nonnegative.
In particular, we will derive the expression for $s$ in (\ref{params-c-s}) for the case $d=nq$ \ ($r=0$).

Bombieri and Davenport \cite{Bombieri} proved that
$$
\mean\limits_{d\in\mathbb N} S(d)
= \mean\limits_{d\in\mathbb N}\prod\limits_{p|d\atop p>2} {p-1\over p-2} = \Pi_2^{-1}
= 1.51478\ldots \ \  (\mbox{\seqnum{A167864}})
$$
where $\Pi_2$ is the {\it twin prime constant}:
$$
\Pi_2=\prod\limits_{p>2}{p(p-2)\over(p-1)^2} = 0.66016\ldots \ \  (\mbox{\seqnum{A005597}}).
$$
If, instead of averaging $S(d)$ over all natural numbers, we restrict the averaging process to an arithmetic progression $r+nq$,
then the resulting average might become larger or smaller than $\Pi_2^{-1}$---or it might remain unchanged ($=\Pi_2^{-1}$).
A few special cases are possible, depending on the values of $q$ and $r$.

\subsection{Base cases: prime $q$}\label{appendix-primeq}

(i) If $q=2$, we have $S(2n)=S(n)$ for all $n \in\mathbb N$.
Therefore, for even $d$,
$$\mean\limits_{d=2n}S(d)=\mean\limits_{d\in\mathbb N}S(d)=\Pi_2^{-1}.
$$
But this immediately implies that, for odd $d$, we must also have
$$\mean\limits_{d=2n+1}S(d)=\mean\limits_{d\in\mathbb N}S(d)=\Pi_2^{-1}.
$$

\rmpara{(ii)}
Suppose $q=p'$ is an odd prime and $\gcd(p',r)=1$. Let
\begin{equation*}
\begin{split}
\xi &\,= \xi(p') = \,\mean\limits_{d=np'} S(d) \\
\eta &\,= \eta(p') = \mean\limits_{d=r+np'} S(d).
\end{split}
\end{equation*}
Observe that every term in progression $d=np'$ has the prime factor $p'$;
so the factor ${p'-1\over p'-2}$ is {\it always present} in the corresponding product (\ref{defSd}).
On the other hand, terms in progression $d=r+np'$ never have the prime factor $p'$;
so the factor ${p'-1\over p'-2}$ is {\it never present} in the corresponding product (\ref{defSd}).
All other odd prime factors $p$ (resp., ${p-1\over p-2}$), $p\ne p'$, are {\it sometimes} present
in both progressions (resp., products), with probability $1/p$.
Therefore, we require that
\begin{equation}\label{xi-eta-S}
\xi=\eta\cdot{p'-1\over p'-2}.
\end{equation}
The factor ${p'-1\over p'-2}$ in (\ref{xi-eta-S}) is the same as in the definition of $S(d)$.

Note that there are $p'-1$ residue classes mod $p'$ for which the mean values of $S(d)$ are $\eta$,
and only one ``special'' residue class (0 mod $p'$) for which the mean value of $S(d)$ is $\xi$.
But the arithmetic average of all these $p'$ mean values is equal to $\mean\limits_{d\in\mathbb N} S(d)$, which is $\Pi_2^{-1}$:
\begin{equation}\label{weighted-avg}
{1\over p'}\cdot\xi + {p'-1\over p'}\cdot \eta = \Pi_2^{-1}.
\end{equation}
Solving equations (\ref{xi-eta-S}) and (\ref{weighted-avg}) together, we find
\begin{align}\label{xi-formula}
 \xi &= \mean\limits_{d=np'} \,S(d)      = {p'\over p'-1} \cdot \Pi_2^{-1}, \\
\eta &=\!\mean\limits_{d=r+np'} \!S(d) = {p'(p'-2)\over(p'-1)^2} \cdot \Pi_2^{-1}. \label{eta-formula}
\end{align}
{\it Example}: for $q=3$, we have
$$\mean\limits_{d=3n} S(d)= \displaystyle{3\over2}\Pi_2^{-1},
\qquad
\mean\limits_{d=3n+1} S(d)= \mean\limits_{d=3n+2} S(d)= \displaystyle{3\over4}\Pi_2^{-1}.
$$
Formulas (\ref{xi-formula}) and (\ref{eta-formula}) encode the following

\bfpara{Factors Principle.} For progressions $d=np'$ where each term is guaranteed to have a factor $p'$,
the average product $S(d)$ equals  $\Pi_2^{-1}$ multiplied by ${p'\over p'-1}>1$.
For progressions $d=r+np'$ where each term {\it never} has the factor $p'$,
the average is  $\Pi_2^{-1}$ multiplied by ${p'(p'-2)\over(p'-1)^2}<1$.

\subsection{Remaining cases: composite $q$}

The above argument considered progressions with one ``special'' prime factor $p'$.
A straightforward modification of the argument allows us to extend the consideration from progressions
with $k$ ``special'' factors to those with $k+1$ ``special'' factors.
Therefore, we can apply the ``factors principle'' to composite values of $q$.

\newpage

\rmpara{(A)}
Let $q$ be a power of 2. (This case is similar to (i) in sect.\,\ref{appendix-primeq}.)
In this case, for every odd prime $p$, terms of our arithmetic progression $r+nq$ span {\it all residue classes} modulo $p$.
So here the average value of $S(d)$ remains unchanged:
$$\mean\limits_{d=r+qn} S(d)=\Pi_2^{-1}.
$$
{\it Example}: for $q=4$, we have
$\mean\limits_{d=4n} S(d)=\mean\limits_{d=4n+1} S(d)= \mean\limits_{d=4n+2} S(d)=\mean\limits_{d=4n+3} S(d)= \Pi_2^{-1}$.

\rmpara{(B)}
Suppose $q\ne2^m$ is composite and $\gcd(q,r)=1$.
Here all terms of our arithmetic progression $r+nq$ are not divisible by prime factors of $q$,
while divisibility by all other primes is neither ensured nor precluded. We have
$$\mean\limits_{d=r+nq} S(d)=\Pi_2^{-1}\cdot \prod\limits_{p|q\atop p>2}{p(p-2)\over(p-1)^2}.
$$
{\it Example}:  for $q=15$ and $r=1$, we have
$\mean\limits_{d=15n+1} S(d)=\displaystyle{3\over4}\cdot{15\over16}\cdot \Pi_2^{-1}=\displaystyle{45\over64} \Pi_2^{-1}$.

\rmpara{(C)}
Suppose $q\ne2^m$ is composite and $r=0$.
All terms $nq$ in our arithmetic progression are divisible by prime factors of $q$,
while divisibility by all other primes is neither ensured nor precluded. We have
$$\mean\limits_{d=nq} S(d)=\Pi_2^{-1}\cdot \prod\limits_{p|q\atop p>2}{p\over p-1}.
$$
{\it Example}: for $q=15$ and $r=0$, we have
$\mean\limits_{d=15n} S(d)=\displaystyle{3\over2}\cdot{5\over4}\cdot \Pi_2^{-1}=\displaystyle{15\over8} \Pi_2^{-1}$.

\rmpara{(D)}
Suppose $q\ne2^m$ is composite, $r\ne0$, $\gcd(q,r)>1$.
(This case is, in a sense, a combination of the above cases.)
Terms of our arithmetic progression $r+nq$ are divisible by common prime factors of $q$ and $r$, but not by other prime factors of $q$;
divisibility by all other primes is neither ensured nor precluded. We have
\begin{equation}\label{meanSd-general}
\mean\limits_{d=r+nq} S(d)=\Pi_2^{-1}
\cdot \prod\limits_{P|q,\,P|r\atop P>2}\!{P\over P-1}\,
\cdot \prod\limits_{p|q,\,p\,\nmid\, r\atop p>2} \!{p(p-2)\over (p-1)^2}.
\end{equation}
{\it Example}: for $q=30$ and $r=3$, we have
$\mean\limits_{d=30n+3} S(d)=\displaystyle{3\over2}\cdot{15\over16}\cdot \Pi_2^{-1}=\displaystyle{45\over32} \Pi_2^{-1}$.

\begin{remark}
Formula (\ref{meanSd-general}) is in fact applicable to an arbitrary pair $(q,r)$;
it subsumes all preceding cases. Formulas like this can be derived not only in our setup related to
the twin prime constant $\Pi_2$ but also for similar average products related to the
Hardy--Littlewood $k$-tuple constants, as well as for other average products, e.g.,
$\mean\limits_{d=r+nq \atop n\in{\mathbb N}} \prod\limits_{p|d\atop p>m}{p+a\over p+b}$.
\end{remark}

\newpage

\section{Acknowledgments}
We are grateful to the anonymous referee and to all contributors and editors
of the websites {\tt MersenneForum.org},  {\tt OEIS.org}, and {\tt PrimePuzzles.net}.
Special thanks to Martin Raab, whose extensive computations provided new examples of
exceptionally large gaps between primes in residue classes \cite{Raab2020}--\cite{Raab2020may}.
We gratefully acknowledge the valuable input and encouragement of the late Prof.~Thomas R.~Nicely,
whose emails we retain to this day.

\bigskip \hrule \bigskip
\noindent {\it 2010 Mathematics Subject Classification}:  Primary 11N05; secondary 11N56, 11N64.
\par
\noindent {\it Keywords}:
arithmetic progression, Brun constant, Cram\'er conjecture,
Gumbel distribution, prime gap, residue class, Shanks conjecture.

\bigskip \hrule \bigskip

\noindent (Concerned with sequences
\seqnum{A000010},
\seqnum{A001620},
\seqnum{A005597},
\seqnum{A014320},
\seqnum{A065421},
\seqnum{A167864},
\seqnum{A194098},
\seqnum{A268925},
\seqnum{A268928},
\seqnum{A330853},
\seqnum{A330854},
\seqnum{A330855},
\seqnum{A334543},
\seqnum{A334544},
\seqnum{A334545},
\seqnum{A335366}, and
\seqnum{A335367}.)

\bigskip
\hrule
\bigskip

\vspace*{+.1in}
\noindent
Published in {\it Journal of Integer Sequences}, {\bf 23} (2020), Article 20.9.3.


\noindent

\end{document}